*Research Article*

# A Generalized Definition of the Fractional Derivative with Applications

M. Abu-Shady[1] and Mohammed K. A. Kaabar[2,3]

[1]*Department of Mathematics and Computer Science, Faculty of Science, Menoufia University, Shibin Al Kawm, Egypt*
[2]*Gofa Camp, Near Gofa Industrial College and German Adebabay, Nifas Silk-Lafto, Addis Ababa 26649, Ethiopia*
[3]*Institute of Mathematical Sciences, Faculty of Science, University of Malaya, Kuala Lumpur 50603, Malaysia*

Correspondence should be addressed to Mohammed K. A. Kaabar; mohammed.kaabar@wsu.edu

A generalized fractional derivative (GFD) definition is proposed in this work. For a differentiable function expanded by a Taylor series, we show that $D^\alpha D^\beta f(t) = D^{\alpha+\beta} f(t); 0 < \alpha \leq 1; 0 < \beta \leq 1$. GFD is applied for some functions to investigate that the GFD coincides with the results from Caputo and Riemann–Liouville fractional derivatives. The solutions of the Riccati fractional differential equation are obtained via the GFD. A comparison with the Bernstein polynomial method (BPM), enhanced homotopy perturbation method (EHPM), and conformable derivative (CD) is also discussed. Our results show that the proposed definition gives a much better accuracy than the well-known definition of the conformable derivative. Therefore, GFD has advantages in comparison with other related definitions. This work provides a new path for a simple tool for obtaining analytical solutions of many problems in the context of fractional calculus.

## 1. Introduction

Fractional calculus theory is a natural extension of the ordinary derivative which has become an attractive topic of research due to its applications in various fields of science and engineering. The integral inequalities in fractional models play an important role in different fields. Massive attention on the advantages of integral inequalities has been paid for considering economics [1], continuum and statistical mechanics [2], solid mechanics [3], electrochemistry [4], biology [5], and acoustics [6]. Fractional-order derivatives of a given function involve the entire function history where the following state of a fractional-order system is not only dependent on its current state but also all its historical states [7, 8]. Nonlocality plays a very important role in several fractional derivative models [9, 10]. Many studies deal with the discrete versions of this fractional calculus by employing the theory of time scales such as [11, 12]. In the literature, some definitions have been introduced such as Riemann–Liouville, Caputo, Jumarie, Hadamard, and Weyl, but all of these definitions have their advantages and disadvantages. The most commonly used definition is Riemann–Liouville which is defined as follows [13].

For $\alpha \in [n-1, n)$, the $\alpha$-derivative of $f(t)$ is

$$D^{RL} f(t) = \frac{1}{\Gamma(n-\alpha)} \frac{d^n}{dx^n} \int_a^\infty \frac{f(x)}{(t-x)^{(\alpha-n+1)}} dx. \quad (1)$$

The Caputo definition is defined as follows.
For $\alpha \in [n-1, n)$, the $\alpha$-derivative of $f(t)$ is

$$D^C f(t) = \frac{1}{\Gamma(n-\alpha)} \int_a^\infty \frac{d^n f(x)/dx^n}{(t-x)^{(\alpha-n+1)}} dx. \quad (2)$$

All definitions including the above (1) and (2) satisfy the linear property of fractional derivatives. These fractional derivatives have several advantages, but they are not suitable for all cases. On the one hand, in the Riemann–Liouville type, when the fractional differential equations are used to describe real-world processes, the Riemann–Liouville derivative has some drawbacks. The Riemann–Liouville derivative of a constant is not zero. Additionally, if an arbitrary function is a constant at the origin, its fractional derivation has a singularity at the origin for instant exponential and Mittag-Leffler functions. Due to these drawbacks, the applicability range of Riemann–Liouville fractional derivatives is limited. For differentiability, Caputo derivative requires higher regularity conditions: to calculate the fractional derivative of a function in the Caputo type, we should first obtain its derivative. Caputo derivatives are defined only for differentiable functions, while the functions that do not have first-order derivative may have fractional derivatives of all orders less than one in the Riemann–Liouville sense (see [13]).

In [14], a new well-behaved simple fractional derivative, named conformable derivative, was defined by relying only on the basic limit definition of the derivative. The conformable derivative satisfies some important properties that cannot be satisfied in Riemann–Liouville and Caputo definitions. However, in [15], the author proved that the conformable definition in [14] cannot provide good results in comparison with the Caputo definition for some functions.

This work aims to provide a new generalized definition of the fractional derivative that has advantages in comparison with other previous definitions in order to obtain simple solutions of fractional differential equations.

The paper is organized as follows: in Section 2, the basic definitions and tools are introduced. In Section 3, some applications are presented. In Section 4, the conclusion is given.

## 2. Basic Definitions and Tools

*Definition 1.* For a function $f: (0, \infty) \longrightarrow R$, the generalized fractional derivative of order $0 < \alpha \leq 1$ of $f(t)$ at $t > 0$ is defined as

$$D^{\text{GFD}} f(t) = \lim_{\varepsilon \longrightarrow 0} \frac{f\left(t + \Gamma(\beta)/\Gamma(\beta - \alpha + 1)\varepsilon t^{1-\alpha}\right) - f(t)}{\varepsilon}; \quad \beta > -1, \beta \in R^+, \tag{3}$$

and the fractional derivative at 0 is defined as $D^{\text{GFD}} f(0) = \lim_{\varepsilon \longrightarrow 0^+} D^{GFD} f(t)$.

**Theorem 1.** *If $f(t)$ is an $\alpha$–differentiable function, then $D^{GFD} f(t) = \Gamma(\beta)/\Gamma(\beta - \alpha + 1)t^{1-\alpha} df(t)/dt; \quad \beta > -1, \beta \in R^+$.*

*Proof.* By using the definition in equation (3), we have

$$D^{\text{GFD}} f(t) = \lim_{\varepsilon \longrightarrow 0} \frac{f\left(t + \Gamma(\beta)/\Gamma(\beta - \alpha + 1)\varepsilon t^{1-\alpha}\right) - f(t)}{\varepsilon};$$
$$\beta > -1, \beta \in R^+, \tag{4}$$

where at $\alpha = \beta = 1$, the classical limit of a derivative function is obtained. Now, let

$$h = \frac{\Gamma(\beta)}{\Gamma(\beta - \alpha + 1)} \varepsilon t^{1-\alpha}, \tag{5}$$

$$\varepsilon = \frac{\Gamma(\beta - \alpha + 1)}{\Gamma(\beta)} h t^{\alpha - 1}. \tag{6}$$

By substituting equation (6) into equation (4), we get

$$D^{\text{GFD}} f(t) = \frac{\Gamma(\beta)}{\Gamma(\beta - \alpha + 1)} t^{1-\alpha} \lim_{h \longrightarrow 0} \frac{f(t+h) - f(t)}{h}. \tag{7}$$

Thus,

$$D^{\text{GFD}} f(t) = \frac{\Gamma(\beta)}{\Gamma(\beta - \alpha + 1)} t^{1-\alpha} \frac{df(t)}{dt}. \tag{8}$$

For a function $f(t) = t^k, k > -1, k \in R^+$, we prove that

$$D^{\text{GFD}} f(t) = \frac{\Gamma(\beta + 1)}{\Gamma(\beta - \alpha + 1)} t^{\beta - \alpha}. \tag{9}$$

By using equation (8), we obtain

$$D^{\text{GFD}} f(t) = \frac{\Gamma(\beta)}{\Gamma(\beta - \alpha + 1)} t^{1-\alpha} k t^{k-1},$$
$$D^{\text{GFD}} f(t) = \frac{k \Gamma(\beta)}{\Gamma(\beta - \alpha + 1)} t^{k - \alpha}. \tag{10}$$

By taking $k = \beta$, we get

$$D^{\text{GFD}} t^\beta = \frac{\beta \Gamma(\beta)}{\Gamma(\beta - \alpha + 1)} t^{\beta - \alpha}, \tag{11}$$

and then

$$D^{\text{GFD}} t^\beta = \frac{\Gamma(\beta + 1)}{\Gamma(\beta - \alpha + 1)} t^{\beta - \alpha}. \tag{12}$$

Equation (12) is compatible with the results of Caputo and Riemann–Liouville derivatives [16]. □

**Theorem 2.** *For a function derivative of $f(t) = t^k$, $k \in R^+$, we obtain $D^\alpha D^\beta t^k = D^{\alpha + \beta} t^k$.*

*Proof.* By using equation (12), we get

$$D^\beta t^k = \frac{\Gamma(k+1)}{\Gamma(k - \beta + 1)} t^{k - \beta},$$

$$D^\alpha D^\beta t^k = \frac{\Gamma(k+1)}{\Gamma(k - \beta + 1)} D^\alpha t^{k - \beta}, \tag{13}$$

$$D^\alpha D^\beta t^k = \frac{\Gamma(k+1)}{\Gamma(k - \beta + 1)} \frac{\Gamma(k - \beta + 1)}{\Gamma(k - \beta - \alpha + 1)} t^{k - \beta - \alpha}.$$

$$\text{L.H.S} = D^\alpha D^\beta t^k = \frac{\Gamma(k+1)}{\Gamma(k - \beta - \alpha + 1)} t^{k - \beta - \alpha}. \tag{14}$$

Also, we have

$$\text{R.H.S} = D^{\alpha + \beta} t^k = \frac{\Gamma(k+1)}{\Gamma(k - \beta - \alpha + 1)} t^{k - \beta - \alpha}. \tag{15}$$

Thus, by (14) and (15), we get

$$D^\alpha D^\beta t^k = D^{\alpha+\beta} t^k. \tag{16}$$

This property is not satisfied in the conformable derivative [14]. □

**Theorem 3.** *For a differentiable function $f(t)$ that expands about a point such as $f(t) = \sum_{k=0}^{\infty} f^k(0)/k! t^k$, we have $D^\alpha D^\beta f(t) = D^{\alpha+\beta} f(t)$.*

*Proof.* The expanded function by Taylor theory is given by $f(t) = \sum_{k=0}^{\infty} f^k(0)/k! t^k$,

$$D^\beta f(t) = \sum_{k=0}^{\infty} \frac{f^k(0)}{k!} D^\beta t^k,$$

$$D^\beta f(t) = \sum_{k=0}^{\infty} \frac{f^k(0)}{k!} \frac{\Gamma(k+1)}{\Gamma(k-\beta+1)} t^{k-\beta},$$

$$D^\alpha D^\beta f(t) = \sum_{k=0}^{\infty} \frac{f^k(0)}{k!} \frac{\Gamma(k+1)}{\Gamma(k-\beta+1)} D^\alpha t^{k-\beta},$$

$$D^\alpha D^\beta f(t) = \sum_{k=0}^{\infty} \frac{f^k(0)}{k!} \frac{\Gamma(k+1)}{\Gamma(k-\beta+1)} \frac{\Gamma(k-\beta+1)}{\Gamma(k-\beta-\alpha+1)} t^{k-\beta-\alpha}, \tag{17}$$

$$\text{L.H.S} = D^\alpha D^\beta f(t) = \sum_{k=0}^{\infty} \frac{f^k(0)}{k!} \frac{\Gamma(k+1)}{\Gamma(k-\beta-\alpha+1)} t^{k-\beta-\alpha}, \tag{18}$$

$$\text{R.H.S} = D^{\alpha+\beta} f(t) = \sum_{k=0}^{\infty} \frac{f^k(0)}{k!} D^{\alpha+\beta} t^k, \tag{19}$$

$$\text{R.H.S} = D^{\alpha+\beta} f(t) = \sum_{k=0}^{\infty} \frac{f^k(0)}{k!} \frac{\Gamma(k+1)}{\Gamma(k-\beta-\alpha+1)} t^{k-\beta-\alpha}. \tag{20}$$

Thus, by equations (18) and (20), we have

$$D^\alpha D^\beta f(t) = D^{\alpha+\beta} f(t). \tag{21}$$

This property is not satisfied in the conformable derivative [14]. □

**Theorem 4.** *Let $\alpha \in (0,1]$ and $f, g$ be $\alpha$–differentiable functions; then,*

$$(i)\ D^{\text{GFD}}(fg) = f D^{\text{GFD}}(g) + g D^{\text{GFD}}(f), \tag{22}$$

$$(ii)\ D^{\text{GFD}}\left(\frac{f}{g}\right) = \frac{g D^{\text{GFD}}(f) - f D^{\text{GFD}}(g)}{g^2}. \tag{23}$$

*Proof.* By using equation (8), we have

$$\text{L.H.S} = D^{\text{GFD}}(fg), \tag{24}$$

$$= \frac{\Gamma(\beta)}{\Gamma(\beta-\alpha+1)} t^{1-\alpha} \frac{d(fg)}{dt}, \tag{25}$$

$$= \frac{\Gamma(\beta)}{\Gamma(\beta-\alpha+1)} t^{1-\alpha} \left[ f \frac{dg}{dt} + g \frac{df}{dt} \right], \tag{26}$$

$$= f \frac{\Gamma(\beta)}{\Gamma(\beta-\alpha+1)} t^{1-\alpha} \frac{dg}{dt} + g \frac{\Gamma(\beta)}{\Gamma(\beta-\alpha+1)} t^{1-\alpha} \frac{df}{dt}, \tag{27}$$

$$= F D^{\text{GFD}}(g) + g D^{\text{GFD}}(f) = \text{R.H.S}. \tag{28}$$

This proves (i).
Now, to prove (ii), we use equation (8) as follows:

$$\begin{aligned}
\text{L.H.S} &= D^{\text{GFD}}\left(\frac{f}{g}\right) \\
&= \frac{\Gamma(\beta)}{\Gamma(\beta-\alpha+1)} t^{1-\alpha} \frac{d}{dt}\left(\frac{f}{g}\right), \\
&= \frac{\Gamma(\beta)}{\Gamma(\beta-\alpha+1)} t^{1-\alpha} \left[\frac{g\, df/dt - f\, dg/dt}{g^2}\right], \\
&= \frac{g\left[\Gamma(\beta)/\Gamma(\beta-\alpha+1) t^{1-\alpha} df/dt\right] - f\left[\Gamma(\beta)/\Gamma(\beta-\alpha+1) t^{1-\alpha} dg/dt\right]}{g^2}, \\
&= \frac{g D^{\text{GFD}}(f) - f D^{\text{GFD}}(g)}{g^2} = \text{R.H.S}.
\end{aligned} \tag{29}$$

Rules (*i*) and (*ii*) are not satisfied in the Caputo and Riemann–Liouville definitions. □

**Theorem 5.** *(Rolle's theorem for the generalized fractional differential function). Let a > 0 and f: [a, b] ⟶ R be a given function that satisfies the following:*

(i) *f is continuous on [a, b]*

(ii) *f is α–differentiable for some α ∈ (0, 1]*

(iii) $f(a) = f(b)$

*Then, there exists $c \in [a, b]$ such that $f^{(\alpha)}(c) = 0$.*

*Proof.* Since $f$ is continuous on $[a, b]$ and $f(a) = f(b)$, there is $c \in (a, b)$, which is a point of local extrema, and $c$ is assumed to be a point of local minimum. So, we have

$$D^{\text{GFD}} f(c^+) = \lim_{\varepsilon \to 0^+} \frac{f\left(c + \Gamma(\beta)/\Gamma(\beta - \alpha + 1)\varepsilon c^{1-\alpha}\right) - f(c)}{\varepsilon}; \quad \beta > -1, \beta \in R^+,$$

$$D^{\text{GFD}} f(c^-) = \lim_{\varepsilon \to 0^-} \frac{f\left(c + \Gamma(\beta)/\Gamma(\beta - \alpha + 1)\varepsilon c^{1-\alpha}\right) - f(c)}{\varepsilon}; \quad \beta > -1, \beta \in R^+.$$

(30)

However, $D^{\text{GFD}} f(c^+)$ and $D^{\text{GFD}} f(c^-)$ have opposite signs. Hence, $D^{\text{GFD}} f(c) = 0$. □

**Theorem 6.** *(mean value theorem for the generalized fractional differential function). Let a > 0 and f: [a, b] ⟶ R be a given function that satisfies the following:*

(i) *f is continuous on [a, b]*

(ii) *f is α–differentiable for some α ∈ (0, 1)*

*Then, there exists $c \in [a, b]$ such that*

$$D^{\text{GFD}} f(c) = \left[\frac{f(b) - f(a)}{h(b^\alpha - a^\alpha)}\right]. \quad (31)$$

*Proof.* Consider a function such as in [25].

$$g(t) = f(t) - f(a) - \left[\frac{f(b) - f(a)}{h(b^\alpha - a^\alpha)}\right](ht^\alpha - ha^\alpha), \quad (32)$$

where $h = 1/\Gamma(\alpha)$.

$$D^{\text{GFD}} g(t) = D^{\text{GFD}} f(t) - D^{\text{GFD}} f(a)$$
$$- \left[\frac{f(b) - f(a)}{h(b^\alpha - a^\alpha)}\right]\left(hD^{\text{GFD}} t^\alpha - hD^{\text{GFD}} a^\alpha\right). \quad (33)$$

By using equation (8), we get

$$D^{\text{GFD}} g(t) = D^{\text{GFD}} f(t) - \left[\frac{f(b) - f(a)}{h(b^\alpha - a^\alpha)}\right], \quad (34)$$

at $c \in [a, b]$.

$$D^{\text{GFD}} g(c) = D^{\text{GFD}} f(c) - \left[\frac{f(b) - f(a)}{h(b^\alpha - a^\alpha)}\right], \quad (35)$$

and the auxiliary function $g(c)$ satisfies all conditions of Theorem 5. Therefore, there exists $c \in [a, b]$ such that $D^{\text{GFD}} g(c) = 0$. Then, we have

$$D^{\text{GFD}} f(c) = \left[\frac{f(b) - f(a)}{h(b^\alpha - a^\alpha)}\right]. \quad (36)$$
□

**Definition 2.** $I_\alpha^a(f)(t) = I_1^0(t^{\alpha-1} f(x)) = \Gamma(\beta - \alpha + 1)/\Gamma(\beta) \int_0^t f(x)/x^{1-\alpha} dx$ and $\alpha \in (0, 1)$.

**Theorem 7.** $D^\alpha I_\alpha(f)(t) = f(t)$ *for $t \geq 0$ where $f$ is any continuous function in the domain.*

*Proof.* Since $f$ is continuous, $I_\alpha^a(f)(t)$ is differentiable. Hence,

$$D^\alpha I_\alpha(f)(t) = \frac{\Gamma(\beta)}{\Gamma(\beta - \alpha + 1)} t^{1-\alpha} \frac{d}{dt} I_\alpha(f)(t),$$

$$D^\alpha I_\alpha(f)(t) = \frac{\Gamma(\beta)}{\Gamma(\beta - \alpha + 1)} t^{1-\alpha} \frac{d}{dt} \frac{\Gamma(\beta - \alpha + 1)}{\Gamma(\beta)} \int_0^\infty \frac{f(x)}{x^{1-\alpha}} dx,$$

$$D^\alpha I_\alpha(f)(t) = t^{1-\alpha} \frac{d}{dt} \int_0^t \frac{f(x)}{x^{1-\alpha}} dx,$$

$$D^\alpha I_\alpha(f)(t) = t^{1-\alpha} \frac{f(t)}{t^{1-\alpha}},$$

$$D^\alpha I_\alpha(f)(t) = f(t).$$

(37)
□

## 3. Applications

### 3.1. Fractional Derivative of the Exponential Function $f(t) = e^{\lambda t}, \lambda \in c$.

$$e^{\lambda t} = \sum_{k=0}^{\infty} \frac{\lambda^k}{k!} t^k,$$

$$D^{\text{GFD}} e^{\lambda t} = \sum_{k=0}^{\infty} \frac{\lambda^k}{k!} D^{\text{GFD}} t^k.$$

(38)

From equation (12), we get

$$D^{\text{GFD}} t^k = D^C t^k. \tag{39}$$

Let us now write equation (23) as

$$D^{\text{GFD}} e^{\lambda t} = \sum_{k=0}^{\infty} \frac{\lambda^k}{k!} D^C t^k, \tag{40}$$

$$D^{\text{GFD}} e^{\lambda t} = D^C e^{\lambda t}.$$

### 3.2. Fractional Derivative of Sine and Cosine Functions.

For the sine function, we define $f(t) = \sin \omega t$ as

$$\sin \omega t = \frac{1}{2!} \left( e^{i\omega t} - e^{-i\omega t} \right),$$

$$D^{\text{GFD}} \sin \omega t = \frac{1}{2!} \left( D^{\text{GFD}} e^{i\omega t} - D^{\text{GFD}} e^{-i\omega t} \right). \tag{41}$$

From equation (26), we obtain

$$D^{\text{GFD}} \sin \omega t = \frac{1}{2!} \left( D^C e^{i\omega t} - D^C e^{-i\omega t} \right),$$

$$^{\text{GFD}} D^\alpha \sin \omega t = {}^C D^\alpha \frac{1}{2!} \left( e^{i\omega t} - e^{-i\omega t} \right), \tag{42}$$

$$^{\text{GFD}} D^\alpha \sin \omega t = {}^C D^\alpha \sin \omega t.$$

Similarly, we can prove the following for $f(t) = \cos \omega t$:

$$^{\text{GFD}} D^\alpha \cos \omega t = {}^C D^\alpha \cos \omega t. \tag{43}$$

Let us now solve some fractional differential equations in the sense of GFD.

*Example 1.* $D^{1/2} y(x) = e^{kx}, \quad y(0) = 0.$

*Solution 1.* Let us find the solution of the above example where $e^{kx} = \sum_{n=0}^{\infty} k^n/n! x^n$.

By using equation (8), we obtain

$$D^{1/2} y(x) = e^{kx},$$

$$\frac{\Gamma(\beta)}{\Gamma(\beta + 1/2)} x^{1/2} \frac{dy(x)}{dx} = \sum_{n=0}^{\infty} \frac{k^n}{k!} x^n,$$

$$\frac{dy(x)}{dx} = \frac{\Gamma(\beta + 1/2)}{\Gamma(\beta)} \sum_{n=0}^{\infty} \frac{k^n}{n!} x^{n-1/2},$$

$$\int dy(x) = \frac{\Gamma(\beta + 1/2)}{\Gamma(\beta)} \sum_{n=0}^{\infty} \frac{k^n}{n!} \int x^{n-1/2} dx,$$

$$y(x) = \frac{\Gamma(\beta + 1/2)}{\Gamma(\beta)} \sum_{n=0}^{\infty} \frac{k^n}{n!} \frac{x^{n+1/2}}{n+1/2} + c,$$

$$y(x) = \sum_{n=0}^{\infty} \frac{k^n}{n!} \frac{\Gamma(\beta + 1/2)}{\Gamma(\beta)} \frac{x^{n+1/2}}{n+1/2} + c. \tag{44}$$

By taking $\beta = n + 1/2$, we have

$$y(x) = \sum_{n=0}^{\infty} \frac{k^n}{n!} \frac{\Gamma(n+1)}{(n+1/2)\Gamma(n+1/2)} x^{n+1/2} + c,$$

$$y(x) = \sum_{n=0}^{\infty} \frac{k^n}{\Gamma(n+3/2)} x^{n+1/2} + c, \tag{45}$$

since $y(0) = 0$.

$$y(x) = \sum_{n=0}^{\infty} \frac{k^n}{\Gamma(n+3/2)} x^{n+1/2}. \tag{46}$$

This solution is consistent with the Caputo solution.

*Example 2.* $D^{1/2} y(x) = x^2 \sin x, \quad y(0) = 0.$

*Solution 2.* Let us find the solution of the above example where $\sin x = \sum_{n=0}^{\infty} x^{2n+1}/(2n+1)!$.

By applying equation (1), we get

$$\frac{\Gamma(\beta)}{\Gamma(\beta + 1/2)} x^{1/2} \frac{dy(x)}{dx} = \sum_{n=0}^{\infty} \frac{x^{2n+3}}{(2n+1)!},$$

$$\frac{dy(x)}{dx} = \frac{\Gamma(\beta + 1/2)}{\Gamma(\beta)} \sum_{n=0}^{\infty} \frac{x^{2n+5/2}}{(2n+1)!},$$

$$\int dy = \frac{\Gamma(\beta + 1/2)}{\Gamma(\beta)} \sum_{n=0}^{\infty} \int \frac{x^{2n+5/2}}{(2n+1)!} dx,$$

$$y(x) = \sum_{n=0}^{\infty} \frac{\Gamma(\beta + 1/2)}{\Gamma(\beta)} \frac{x^{2n+7/2}}{(2n+7/2)(2n+1)!} + c. \tag{47}$$

By taking $\beta = 2n + 7/2$, we get

$$y(x) = \sum_{n=0}^{\infty} \frac{\Gamma(2n+4)}{\Gamma(2n+7/2)} \frac{x^{2n+7/2}}{(2n+7/2)(2n+1)!} + c,$$

$$y(x) = \sum_{n=0}^{\infty} \frac{(2n+3)!}{(2n+7/2)\Gamma(2n+7/2)} \frac{x^{2n+7/2}}{(2n+1)!} + c, \tag{48}$$

$$y(x) = \sum_{n=0}^{\infty} \frac{(2n+3)(2n+2)}{\Gamma(2n+9/2)} x^{2n+7/2} + c,$$

since $y(0) = 0$.

$$y(x) = \sum_{n=0}^{\infty} \frac{(2n+3)(2n+2)}{\Gamma(2n+9/2)} x^{2n+7/2}. \tag{49}$$

This solution is consistent with the Caputo solution.

*Example 3.* $D^{1/2} y(x) + y(x) = x^2 + 2/\Gamma(2.5) x^{3/2}.$

*Solution 3.* By applying equation (8), we obtain

$$\frac{\Gamma(\beta)}{\Gamma(\beta + 1/2)} x^{1/2} \frac{dy}{dx} + y(x) = x^2 + \frac{2}{\Gamma(2.5)} x^{3/2}. \tag{50}$$

Table 1: Comparison of the results of the GFD with other works at $\alpha = 3/4$.

| t | Present work | BPM [17] | EHPM [18] | IABMM [18] | CD [14] |
|---|---|---|---|---|---|
| 0 | 0 | 0 | 0 | 0 | 0 |
| 0.2 | 0.31439 | 0.30996891 | 0.3214 | 0.3117 | 0.37889 |
| 0.4 | 0.49848 | 0.48162749 | 0.5077 | 0.4855 | 0.58539 |
| 0.6 | 0.63022 | 0.59777979 | 0.6259 | 0.6045 | 0.72064 |
| 0.8 | 0.72609 | 0.67884745 | 0.7028 | 0.6880 | 0.81029 |
| 1.0 | 0.79618 | 0.73684181 | 0.7542 | 0.7478 | 0.87006 |

Table 2: Comparison of the results of the GFD with other works at $\alpha = 9/10$.

| t | Present work at $\alpha = 9/10$ | BPM [17] | MHPM [18] | IABMM [18] | CD [14] |
|---|---|---|---|---|---|
| 0 | 0 | 0 | 0 | 0 | 0 |
| 0.2 | 0.23952 | 0.23878798 | 0.2391 | 0.2393 | 0.25526 |
| 0.4 | 0.42667 | 0.42258214 | 0.4229 | 0.4234 | 0.45191 |
| 0.6 | 0.57607 | 0.56617082 | 0.5653 | 0.5679 | 0.60539 |
| 0.8 | 0.69138 | 0.67462642 | 0.6740 | 0.6774 | 0.72063 |
| 1.0 | 0.7778 | 0.75460256 | 0.7569 | 0.7584 | 0.80445 |

Table 3: Comparison of the results of the GFD with other works at $\alpha = 9/10$ for equation (43).

| t | Present work | BPM [17] | FTBM [18] | IABMM [18] | CD [14] |
|---|---|---|---|---|---|
| 0 | 0 | 0 | 0 | 0 | 0 |
| 0.2 | 0.30718 | 0.31488815 | 0.31485423 | — | 0.33295 |
| 0.4 | 0.67131 | 0.69756771 | 0.69751826 | — | 0.73105 |
| 0.6 | 1.0666 | 1.10789047 | 0.90364539 | — | 1.1561 |
| 0.8 | 1.4397 | 1.47772823 | 1.47768008 | — | 1.5422 |
| 1.0 | 1.7485 | 1.76542008 | 1.76525852 | 1.7356 | 1.8457 |

The following is a nonlinear differential equation of first order in which we can obtain its solution with the help of Mathematica package.

$$y(x) = \frac{\begin{array}{c}9A^4\sqrt{\pi} - 6A^3(2 + 3\sqrt{\pi}\sqrt{x})\sqrt{\pi}\sqrt{x} - 12A(2 + \sqrt{\pi}\sqrt{x})x + 16x^{3/2} \\ + 6\sqrt{\pi}x^2 + 6A^2(4\sqrt{x} + 3\sqrt{\pi}x)\end{array}}{6\sqrt{\pi}} + c_1 e^{-2\sqrt{x}/A}. \qquad (51)$$

To determine $A = \Gamma(\beta)/(\Gamma(\beta + 1/2))$, we take $\alpha = \beta = 1/2$ as in [19].

*Example 4.* $d/dx\{(1 - \sqrt{x})(y(x) + 1)\} + \lambda D^{1/2} y(x) = 0$.

*Solution 4.* By applying equation (8), we obtain

$$\frac{d}{dx}\{(1 - \sqrt{x})(y(x) + 1)\} + \lambda \frac{\Gamma(\beta)}{\Gamma(\beta + 1/2)} x^{1/2} \frac{dy}{dx} = 0. \qquad (52)$$

The following is a nonlinear differential equation of first order in which its solution can be obtained using Mathematica package as mentioned in our previous example.

*Example 5.* Consider the fractional Riccati differential equation [17]:

$$D^\alpha y(x) + y^2(x) = 1, y(0) = 0, 0 < \alpha \leq 1. \qquad (53)$$

*Solution 5.* By applying equation (8), we obtain

$$\frac{\Gamma(\beta)}{\Gamma(\beta - \alpha + 1)} x^{1-\alpha} \frac{dy}{dx} + y^2(x) = 1, y(0) = 0, 0 < \alpha \leq 1. \qquad (54)$$

To solve this equation at $\alpha = 3/4$ and $\alpha = 9/10$, the package of Mathematica has been used to obtain the following:

$$y(x) = \frac{-1 + e^{8x^{3/4}/sA}}{1 + e^{8x^{3/4}/sA}}, \qquad (55)$$

where $A = \Gamma(\beta)/\Gamma(\beta + 1/4)$ and $\beta = \alpha = 3/4$ as in [28].

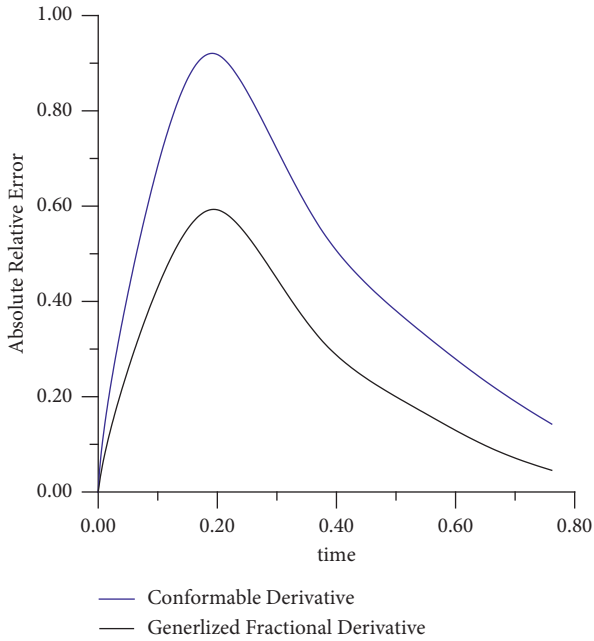

Figure 1: The absolute relative error is plotted for the Riccati fractional differential equation for the conformable derivative and GFD at $\alpha = 0.75$.

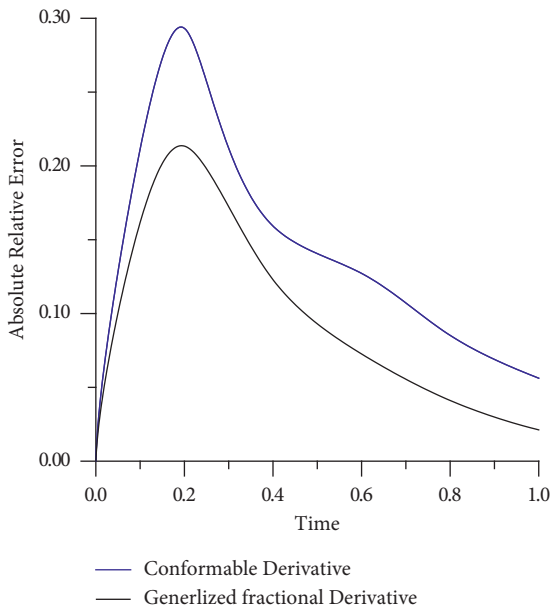

Figure 2: The absolute relative error is plotted for the Riccati fractional differential equation for the conformable derivative and GFD at $\alpha = 0.90$.

$$y(x) = \frac{-1 + e^{20x^{9/10}/9A}}{1 + e^{20x^{9/10}/9A}}, \quad (56)$$

where $A = \Gamma(\beta)/\Gamma(\beta + 1/4)$ and $\beta = \alpha = 9/10$ as in [19].

*Example 6.* Consider the following Riccati fractional differential equation [17]:

$$D^\alpha y(x) = 2y(x) - y^2(x) + 1, \; y(0) = 0, 0 < \alpha \le 1. \quad (57)$$

*Solution 6.* By applying equation (8), we obtain

$$\frac{\Gamma(\beta)}{\Gamma(\beta - \alpha + 1)} x^{1-\alpha} \frac{dy}{dx} = 2y(x) - y^2(x) + 1, \; y(0) = 0, 0 < \alpha \le 1. \quad (58)$$

To solve this equation at $\alpha = 9/10$, the package of Mathematica has been used to obtain the following:

$$y(x) = -\left(\frac{-1 - \sqrt{2} - e^S + \sqrt{2}e^S}{1 + e^S}\right), \quad (59)$$

where

$$s = \frac{2\sqrt{2}\left(-10x^{9/10} + 9A \ln(1 + \sqrt{2})/\sqrt{2}\right)}{9A}, \quad (60)$$

$A = \Gamma(\beta)/\Gamma(\beta + 1/4)$ and $\beta = \alpha = 9/10$ as in [19].

## 4. Discussion of Results

In this section, we show some results for the Riccati fractional differential equation in Tables 1–3 for different values of $\alpha$, where parameters are taken as $\beta = \alpha$ [19]. As a result, we have obtained a good accuracy in the present calculations, where $\alpha = 3/4$ is taken in Table 1, and $\alpha = 9/10$ is taken in Tables 2 and 3. By comparing our results from the GFD definition with the Bernstein polynomial method (BPM) [17], enhanced homotopy perturbation method (EHPM) [18], IABMM [12], and conformable derivative (CD) [14], it is noticeable that the present results are in good agreement with BPM, EHPM, and IABMM results. In addition, the conformable derivative [14] has been used to solve the fractional Riccati differential equation. However, the results of the conformable derivative do not coincide with other works and our present results. A similar situation is in Table 2, by taking $\alpha = 9/10$, where the present results are compared with the Bernstein polynomial method (BPM) [17], enhanced homotopy perturbation method (EHPM) [18], IABMM [18], and conformable derivative (CD). The obtained results that have been calculated analytically via the GFD are in good agreement with other methods. However, in comparison with the CD, the present results are better than CD results as suggested in [14]. In Figure 1, the absolute relative error shows that the present results of the Riccati fractional differential equation are exactly obtained at $\alpha = 1$ in [17] by comparing it with $\alpha = 3/4$ using the proposed definition and the conformable one. The figure shows a good accuracy for the results of the proposed definition in comparison with the conformable one. A similar situation is provided in Figure 2 at $\alpha = 9/10$.

...

...

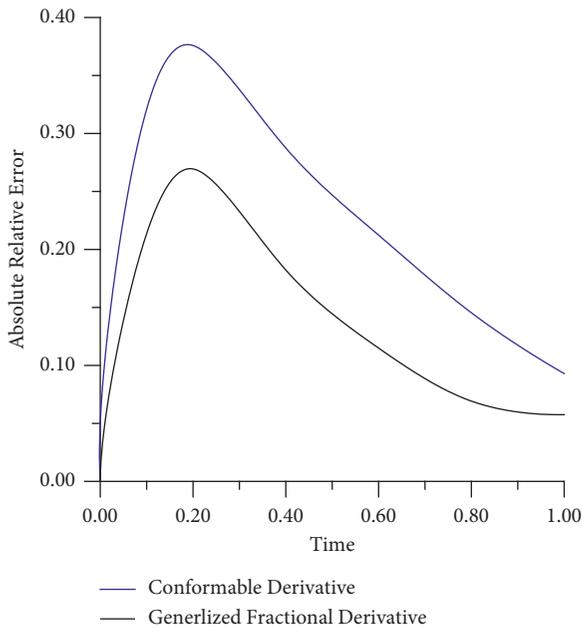

Figure 3: The absolute relative error is plotted for the Riccati fractional differential equation for the conformable derivative and GFD at $\alpha = 0.90$.

In Table 3, one compares the present results obtained from the GFD definition with the Bernstein polynomial method (BPM) [6], fractional Taylor basis method (FTBM) [7], IABMM [8], and conformable derivative (CD) [2]. The results from numerical methods in [6–8] are in agreement with the present results in the context of GFD. In comparison with the conformable derivative results, CD gives more errors than our results that have been obtained in the sense of the GFD. Also, in Figure 3, the results of the GFD give less error in comparison with the conformable derivative definition. Therefore, the present results of the GFD give compatible results with other works.

## 5. Conclusion

In this work, GFD has been suggested to provide more advantages than other classical Caputo and Riemann–Liouville definitions such as the derivative of two functions, the derivative of the quotient of two functions, Rolle's theorem, and the mean value theorem which have been satisfied in the GFD. The present definition satisfies $D^\alpha D^\beta f(t) = D^{\alpha+\beta} f(t)$ for a differentiable function $f(t)$ expanded by Taylor series. The fractional integral is introduced. Compatible results with Caputo and Riemann–Liouville results have been obtained for functions that are given in Sections 3.1 and 3.2. Also, a comparison with the conformable derivative is studied.

Some fractional differential equations can be solved analytically in a simple way with the help of our proposed definition which exactly agrees with the classical Caputo and Riemann–Liouville derivatives' results. In comparison with the conformable derivative, less error has been obtained in our GFD results by calculating the absolute relative error as in Figures 1–3 for the given Riccati fractional differential equation. Also, our results from the GFD definition are compared with the Bernstein polynomial method (BPM), enhanced homotopy perturbation method (EHPM), IABMM, and conformable fractional derivative (CD) [14]. The present results are in good agreement with BPM, EHPM, and IABMM.

We conclude that the present definition gives a new direction for solving fractional differential equations in a simple manner in which the results of the Caputo and Riemann–Liouville definitions are exactly deduced. In addition, GFD has advantages in comparison with the conformable derivative definition.

## Data Availability

No data were used to support this study.

## Conflicts of Interest

The authors declare that they have no conflicts of interest.